\documentclass[12pt,a4paper]{article}
\usepackage[latin2]{inputenc}
\usepackage{graphicx}
\usepackage{amsmath}
\usepackage{amsfonts}
\begin{document}
\begin{center}\Large{\bf{ON THE TOTAL GRAPH OF A FINITE COMMUTATIVE RING
}}\end{center}

\begin{center}\textbf{\large M. H. Shekarriz\footnote{\small E-mail addresses: mhshekarriz@gmail.com (M. H.
Shekarriz),\\ shirdareh@susc.ac.ir (M. H. Shirdareh Haghighi),
sharif@susc.ac.ir (H. Sharif).}, M. H. Shirdareh Haghighi, H.
Sharif}\end{center}

\begin{center}\textit{\bf Department of Mathematics, College of
Science,\\
Shiraz University, 71454,
 Shiraz, Iran}\\
\end{center}

\begin{abstract}{\sl Let $R$ be a finite commutative ring with
$1\ne 0$. In this article, we study the total graph of $R$,
denoted by $\tau (R)$, determine some of its basic
graph-theoretical properties, determine when it is Eulerian, and
find some conditions under which this graph is isomorphic to
$Cay(R,Z(R)\backslash\lbrace 0\rbrace)$. We shall also compute
the domination number of $\tau (R).$}
\end{abstract}
\begin{center} \textit{\small Keywords: Cayley graphs, Total graph of a commutative ring, Zero-divisors} \end{center}

\section{INTRODUCTION}

Let $R$ be a commutative ring with $1\ne 0$, $Z(R)$ its set of
zero-divisors, and $R^{*}$ its set of units. Anderson and Badawi
(2008) introduced \textit{the total graph of} $R$, which we
denote here by $\tau (R)$, as a simple graph with all elements of
$R$ as vertices and for distinct $x,y\in R$, the vertices $x$ and
$y$ are adjacent if and only if $x+y\in Z(R)$. They also studied
some basic properties of $\tau (R)$ such as diameter, girth,
connectivity, etc., mostly without the assumption that $R$ is
finite. Akbari et al. (2009) have shown that if $R$ is a finite
commutative ring and $\tau (R)$ is connected, it is also
Hamiltonian. Maimani et al.  have characterized all finite rings
whose total graphs are planar or toroidal and have also shown
that for every positive integer $g$, there are only finitely many
finite rings (up to isomorphism) whose total graphs have genus
$g$.
\vspace{1mm}\\
\textbf{Remark 1.1.} Anderson and Badawi (2008) use the notation
$T(\Gamma (R))$ instead of $\tau(R)$. We prefer $\tau (R)$ because
$T(\Gamma(R))$ has a different meaning in graph-theoretical
context (see, for example, Behzad (1970)).
\vspace{1mm}\\
\indent The definition of $\tau (R)$ may bring to mind the definition of the
Cayley graph, $Cay(R, Z(R)\backslash\lbrace 0\rbrace )$, which we
denote here by $C(R)$. Hence, hereinafter, $C(R)$ represents the
simple graph with vertex set $R$ and for distinct $x,y\in R$, the
vertices $x$ and $y$ are adjacent if and only if $x-y\in Z(R)$.
Cayley graphs are widely studied in algebraic graph theory, and most of
their graph-theoretical properties are known. For instance, it is known
that $C(R)$ is a $(\vert Z(R)\vert -1)-$regular graph. For a
general reference on algebraic graph theory, see Beineke and Wilson
(2004). Meanwhile, Akhtar et al. (2009) have also studied the unitary
Cayley graph of a finite ring $R$, $Cay(R, R^{*})$,
which is the complement of $C(R)$.

Besides the resemblance between the definitions of $\tau (R)$ and $
C(R)$, they may be quite different in some of their graph-theoretical
properties. In this article, we answer the naturally arising question:
under what conditions on a finite commutative ring $R$, do we have $
\tau (R)\simeq C(R)$?

We need some well-known facts about commutative rings: If $R$ is
an Artinian ring, then either $R$ is local with its maximal ideal
$\frak{M}$, or $R\cong R_{1}\oplus \cdots  \oplus R_{k}$, where
$k\ge 2$ and each $R_{i}$ is local with maximal ideal
$\frak{M}_{i}$; this decomposition is unique up to permutation of
factors, see Atiyah and Macdonald (1969). Moreover, if $R$ is
finite, then every element of $R$ is either a unit or a
zero-divisor, i.e., $\vert Z(R)\vert =\vert R\vert -\vert
R^{*}\vert $. Furthermore, if $R$ is also a local ring with
maximal ideal $\frak{M}$, then $\frak{M}=Z(R)$, and there exists a
prime $p$ such that the characteristic of the residue field
$R/\frak{M}$ is $p$, and $\vert R\vert $, $\vert \frak{M}\vert $,
and $\vert R/\frak{M}\vert $ are all powers of $ p$. Moreover, if
$R \cong R_{1} \oplus \cdots \oplus R_{k}$, then $(u_{1},\ldots
,u_{k})$ is a unit in $R$ if and only if $u_{i} \in R_{i}^{*}$
for each $i=1,\ldots ,k$, and therefore we have $\vert R^{*}\vert
=\vert R_{1}^{*}\vert \times \cdots \times \vert R_{k}^{*}\vert $.

In this article, we denote the residue field $R_{i}/\frak{M}_{i}$
by $F_{i} $, the quotient map by $\pi _{i}:R_{i} \longrightarrow
F_{i}$, and $\vert F_{i}\vert $ by $f_{i}$. As mentioned in
Akhtar et al. (2009), after appropriate permutation of factors,
we may assume that $f_{1}\le \ldots \le f_{k}$.

For the moment, suppose that $R$ is an infinite commutative ring.
Then either $R$ has more than one zero-divisor, which means that
it has infinitely many zero-divisors, or $R$ is an integral
domain, which means that its total graph is the disjoint union of
infinitely many $ K_{1}$'s or $K_{2}$'s. In the former case, the
degree of each vertex in $\tau (R)$ is infinite, which means that
$\tau(R)$ is not locally finite and most of graph-theoretical
properties such as being Eulerian or the domination number, etc.
are meaningless or seem to be hard to determine. Thus,
hereinafter, we assume that all rings are finite commutative with
$1\ne 0$. Meanwhile, we use the notations and definitions of
graph theory from West (2000).

\section{ BASIC PROPERTIES OF $\boldsymbol{\tau (R)}$}

The following two propositions can be easily proved using theorems of
Anderson and Badawi (2008, 2.2, 2.6, 3.3, 3.4 and 3.14):
\vspace{2mm}\\
\textbf{Proposition 2.1.} \textit{Let $R$ be a finite commutative
ring. Then $\tau (R)$ is connected if and only if $Z(R)$ is not
an ideal of $R$, if and only if $R$ is not a local ring.
Moreover, if $R$ is not a local ring, then {\normalfont
diam}$(\tau (R))=2$, and {\normalfont gr}$(\tau (R))=3$ except
when $R \cong \Bbb{Z}_{2} \oplus {\Bbb{Z}}_2,$ where we have
{\normalfont gr}$(R)=4$.} $\Box$
\vspace{2mm}\\
\textbf{Proposition 2.2. }\textit{Let $R$ be a local ring and
$\beta =\vert R/Z(R)\vert $. Then
 \begin{enumerate}
\item[{\rm(a)}]  if $2\in Z(R)$, then $\tau (R)$ is the union of $\beta $
disjoint $K_{\vert Z(R)\vert }$'s;
\item[{\rm(b)}] if $2 \notin Z(R)$, then
$\tau (R)$ is the disjoint union of one copy of $K_{\vert
Z(R)\vert }$ and $(\beta -1)/2$
copies of $K_{\vert Z(R)\vert ,\vert Z(R)\vert }$.
\end{enumerate}
Moreover, {\normalfont gr}$(\tau (R))=3$ if $\vert Z(R)\vert \ge
3$, otherwise {\normalfont gr}$(\tau (R))=\infty.$} $\Box$
\vspace{2mm}\\
\textbf{Remark 2.3.} Anderson and Badawi (2008, 2.6: (3)(b)) have
stated that when $Z(R)$ is an ideal of $R$, then gr$(\tau (R))=4$
if and only if $2\notin Z(R)$ and $\vert Z(R)\vert =2$. But this
case cannot happen since there are only two (finite) rings with
$\vert Z(R)\vert =2$, say $\Bbb{Z}_{4}$ and $\Bbb{Z}_{2}\lbrack
X\rbrack /(X^{2})$, and in both of them $2 \in Z(R)$. (This
Remark is also clear form Lemmas 2.5 and 2.6 below.)

In $\tau (R)$, the degree of $0$ is $\vert Z(R)\vert -1$.
Maimani et al. have proved a lemma similar to the following:
\vspace{2mm}\\
\textbf{Lemma 2.4.} \textit{Let $R$ be a finite
commutative ring and $x$ be a vertex of $\tau (R)$.
Then the degree of $x$ is $\vert Z(R)\vert $
if $x+x=2x \notin Z(R)$, otherwise the degree of $x$
 is $\vert Z(R)\vert -1$. In particular,
$2\in Z(R)$ if and only if $\tau (R)$ is a $(\vert
Z(R)\vert -1)-$regular graph.}
\vspace{2mm}\\
\textbf{Proof.} For every $z \in Z(R)$, there is a unique
$a=z-x\in R$ such that $ x+a=z$. Then $x$ is adjacent to $a$
unless $a=x$. If $2x\in Z(R) $, then since $x$ cannot be adjacent
to itself, deg$x=\vert Z(R)\vert -1$. If $2x\notin Z(R)$, then for
every zero-divisor $z$ we have $z-x\ne x$, and hence deg$x=\vert
Z(R)\vert $. In particular, if $2\in Z(R)$, then for every $x\in
R$, we have deg$x=\vert Z(R)\vert -1$. Therefore, $\tau (R)$ is a
$(\vert Z(R)\vert -1)$-regular graph. Now suppose that $2\notin
Z(R)$. Then deg$1=\vert Z(R)\vert \ne \vert Z(R)\vert -1=$deg0,
i.e., $\tau (R)$ is not a regular graph. $\Box$

We also need the following two lemmas whose proofs are simple and are
left to the reader.
\vspace{2mm}\\
\textbf{Lemma 2.5.} \textit{Let $R$ be a finite commutative ring.
Then $2\in Z(R)$ if and only if $2 {\big \vert} \vert R\vert.$ }
$\Box$
\vspace{2mm}\\
\textbf{Lemma 2.6.} \textit{Let $R$ be a finite commutative ring.
If $\vert R\vert $ is odd, then $\vert R^{*}\vert $ is even.}
$\Box$

Summing up, we have the following theorem:
\vspace{2mm}\\
\textbf{Theorem 2.7.} \textit{Let $R$ be a finite commutative
ring. Then
\begin{enumerate}
\item[{\rm(a)}] if $\vert R\vert$  is even, then $\tau (R)$ is a $(\vert
Z(R)\vert -1)-$regular graph,
\item[{\rm(b)}] if $\vert R\vert $ is odd and $r\in R^{*}$ (respectively, $r \in
Z(R)$), then {\normalfont deg}$r=\vert Z(R)\vert $ (respectively,
{\normalfont deg}$r=\vert Z(R)\vert -1$). $\Box$
\end{enumerate}}

 \noindent \textbf{Remark 2.8.} Note that $1$ and $0$ are not
adjacent in $ \tau (R)$. Therefore, $\tau (R)$ is never a
complete graph. Furthermore, we show that $\tau (R)$ is never an
odd cycle either. To the contrary, suppose that $\tau (R)$ is an
odd cycle. Hence deg1=deg0=2=$\vert Z(R)\vert -1$. Therefore, by
Lemma 2.4 we have $ 2\in Z(R)$, and hence by Lemma 2.5, we must
have $2 {\big \vert} \vert R\vert$, which is a contradiction
since $\vert R\vert $ is odd.

\section{WHEN IS $\boldsymbol{\tau (R)}$ EULERIAN?}

A graph is said to be \textit{Eulerian} if it has a closed trail
containing all edges. It is well-known that a graph is Eulerian if and
only if it is connected and its vertex degrees are all even, see West
(2000).
\vspace{2mm}\\
\textbf{Lemma 3.1.} \textit{Let $R$ be a finite commutative ring.
Then the following three conditions are necessary for $\tau (R)$
to be an Eulerian graph:
\begin{enumerate}
\item[{\rm(a)}] $Z(R)$ is not an ideal of $R$,
\item[{\rm(b)}] $\vert Z(R)\vert$ is an odd integer,
\item[{\rm(c)}] $2\in Z(R)$.
\end{enumerate}}
\noindent
 \textbf{Proof.} Condition (a) is necessary since otherwise by
Proposition 2.1, $ \tau (R)$ is not connected and cannot be
Eulerian. (b) is also necessary since otherwise deg$0=\vert
Z(R)\vert -1$ is odd, and thus $ \tau (R)$ cannot be Eulerian.
Condition (c) is also necessary, since otherwise by Lemma 2.4 we
have deg$1=\vert Z(R)\vert =1+$deg0, and hence $\tau (R)$ cannot
be an Eulerian graph since it has at least one vertex with an odd
degree. $\Box$
\vspace{2mm}\\
\textbf{Remark 3.2.} The three conditions above are also sufficient
because they imply that $\tau (R)$ is a connected regular graph of an
even degree, i.e., $\tau (R)$ is Eulerian. Moreover, while conditions
(a) and (b) hold, condition (c) above could be replaced by

\begin{enumerate}
\item[{\rm(c$'$)}] \textit{ $1+1=0$, i.e., the characteristic of $R$ is 2.}
\end{enumerate}

Condition (c$'$) always implies (c). Conversely, assume necessary
conditions (a) and (b) hold. By Lemma 2.5, we have $2\in Z(R)$
implies $ 2{\big \vert} \vert R\vert $. Meanwhile, by condition
(b), $\vert Z(R)\vert $ is odd and thus $\vert R^{*}\vert $ is
also odd. Therefore, there exists some $u\in R^{*}$ such that
$u=-u$. Then $uu^{-1}=-uu^{-1}$ which means that $1=-1$. Hence
$2=0$.
\vspace{2mm}\\
\textbf{Theorem 3.3.} \textit{Let $R$ be a finite commutative
ring. Then the graph $\tau (R)$ is Eulerian if and only if $R$ is
isomorphic to a direct sum of two or more finite fields of even
orders, i.e., $R\cong \bigoplus_{i=1}^{k}\Bbb{F}_{2^{t_{i}}}$ for
some $k\ge 2$.}
\vspace{2mm}\\
\textbf{Proof.} If $R\cong
\bigoplus_{i=1}^{k}\Bbb{F}_{2^{t_{i}}}$ for some\textit{ }$k\ge
2$,\textit{ }then $\tau (R)$ is Eulerian by Lemma 3.1 and Remark
3.2. To prove the converse, suppose that $R$ is a ring whose
total graph is Eulerian. Consequently, by Lemma 3.1, $R\cong
R_{1} \oplus \cdots  \oplus R_{k}$, $k\ge 2$. Meanwhile, by
Remark 3.2, we must have $1+1=0$ , which means that
$1_{R_{i}}+1_{R_{i}}=0$ for all $i=1,\ldots ,k$. Hence, by Lemma
2.5, it can be inferred that $\vert R_{i} \vert$'s are all even.
If $\vert Z(R_{i})\vert \ge 2$ for some $i=1,\ldots ,k$, then $2
{\big \vert} \vert Z(R_{i})\vert$, and hence $\vert
R_{i}^{*}\vert =\vert R_{i}\vert -\vert Z(R_{i})\vert $ is even.
Therefore, $\vert R^{*}\vert =\vert R_{1}^{*}\vert \times \cdots
\times \vert R_{k}^{*}\vert $ is also even, i.e., $\vert
Z(R)\vert $ is even which is a contradiction because $\tau (R)$
is an Eulerian graph and by Lemma 3.1 (b), $\vert Z(R)\vert $
must be odd. Thus, we must have $\vert Z(R_{i})\vert =1$ for all
$i=1,\ldots ,k$. Therefore, the $R_{i}$'s are all finite fields
of even orders, i.e., $R\cong
\bigoplus_{i=1}^{k}\Bbb{F}_{2^{t_{i}}}$. $\Box$
\vspace{2mm}\\
\textbf{Example 3.4.} $\tau (\Bbb{Z}_{2} \oplus \Bbb{Z}_{2})$ and
$\tau (\Bbb{Z}_{2} \oplus \Bbb{Z}_{2} \oplus \Bbb{Z}_{2})$ are
shown in figures 1 and 2, respectively. They are both Eulerian:
\begin{figure}[h]
\centering
\includegraphics[width=2.13cm,height=2.13cm]{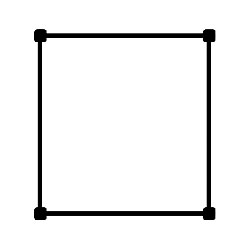}
\end{figure}
\begin{center}\textbf{Figure 1.} $\tau (\Bbb{Z}_{2}\oplus\Bbb{Z}_{2})$.\end{center}
\begin{figure}[h]
\centering
\includegraphics[width=2.66cm,height=2.66cm]{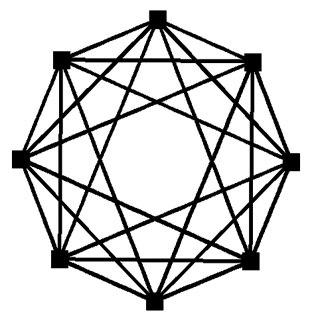}
\end{figure}
\begin{center}\textbf{Figure 2.} $\tau (\Bbb{Z}_{2}\oplus\Bbb{Z}_{2}\oplus\Bbb{Z}_{2})$.
\end{center}

\section{DOMINATION NUMBER}

In a graph $G$, a set $S \subseteq V(G)$ is called a
\textit{dominating set} if every vertex not in $S$ has a neighbor
in $S$. The \textit{domination number} $\gamma (G)$ is the
minimum size of a dominating set in $G$.
\vspace{2mm}\\
\textbf{Theorem 4.1.} \textit{Let $R$ be a finite commutative
ring. We have $\gamma (\tau (R))=f_{1}$, except when $R$ is an
integral domain of an odd order, where $\gamma (\tau
(R))=\frac{f_{1}-1}{2}+1$.}
\vspace{2mm}\\
\textbf{Proof.} Let us treat the local case first. If $R$ is a
local ring and $ \vert R\vert $ is odd, then by Proposition 2.2,
$\tau (R)$ is a disjoint union of one copy of $K_{\vert Z(R)\vert
}$ and $(f_{1}-1)/2 $ copies of $K_{\vert Z(R)\vert ,\vert
Z(R)\vert }$. If besides we have $\vert Z(R)\vert =1$, then every
dominating set must contain the vertex $0$ and at least one
vertex from each copy of $K_{1,1}$, and hence $\gamma (\tau
(R))=\frac{(f_{1}-1)}{2}+1$. If $\vert R\vert $ is odd and $\vert
Z(R)\vert \ne 1$, then every dominating set of a copy of
$K_{\vert Z(R)\vert ,\vert Z(R)\vert }$ must contain at least $2$
vertices, one from each partite set. Therefore, we have $\gamma
(\tau (R))=f_{1}$.

If $R$ is a local ring and $\vert R\vert $ is even, then $\tau
(R)$ is a disjoint union of $f_{1}$ copies of $K_{\vert Z(R)\vert
}$. Hence every dominating set requires at least one vertex from
each copy. Thus $\gamma (\tau (R))=f_{1}$.

Now suppose that $R$ is not a local ring. Then $R\cong R_{1}\oplus
\cdots \oplus R_{k}$; $k\ge 2$. Choose $f_{1}$ elements
$x_{j}=(x_{1j},\ldots ,x_{kj})\in R$, for $j=1,\ldots ,f_{1}$,
such that $\pi _{1}(x_{1j})\ne \pi _{1}(x_{1j^{'}})$ for $j\ne
j'$. Put $S=\lbrace x_{j} \vert j=1,\ldots ,f_{1}\rbrace $, and
suppose that $y=(y_{1},\ldots ,y_{k})$ is an arbitrary element of
$R\setminus S$. Then there is an $i=1,...,f_{1}$ such that $\pi
_{1}(x_{1i})=\pi _{1}(y_{1}).$ If $ f_{1}$ is even, let $j=i$,
else let $j$ be the index of the unique element of $S$, $x_{j}$
say, such that $\pi _{1}(x_{1j})=-\pi _{1}(x_{1i})$. Thus
$y_{1}+x_{1j}\in Z(R_{1})$, and hence $y$ and $ x_{j}$ are
adjacent in $\tau (R)$. Therefore, $S$ is a dominating set and
$\gamma (\tau (R))\le \vert S\vert =f_{1}$.

It remains to prove that $\gamma (\tau (R))\ge f_{1}$. Let $
A=\lbrace a_{i}=(a_{1i},\ldots ,a_{ki})\in R \;\; \vert\;\;
i=1,\ldots,t \mbox{ and } t<f_{1}\rbrace.$ We show that $A$ is
not a dominating set, and this completes the proof. For each
$j=1,\ldots ,k$, since $t<f_{1}\le f_{j}$, there is $b_{j}\in
R_{j}$ such that $\pi _{j}(b_{j})\ne -\pi _{j}(a_{ji})$ for each
$i=1,\ldots ,t$. Put $b=(b_{1},\ldots ,b_{k})$. Then $
b+a_{i}\notin Z(R)$ for each $i=1,\ldots ,t$, which means that
$b$ does not have a neighbor in $A$, and hence $A$ is not a
dominating set. $\Box$
\vspace{2mm}\\
\textbf{Example 4.2.} A dominating set for $\tau (R)$ for $R\cong
\Bbb{Z}_{2}$, $\Bbb{Z}_{4}$, $\Bbb{Z}_{4}\oplus \Bbb{Z}_{3}$, or
$\oplus _{i=1}^{k}\Bbb{Z}_{2}$ is $\lbrace 0_{R},1_{R}\rbrace $. A
dominating set for $\tau (\Bbb{Z}_{45})$ is $\lbrace
\bar{0},\bar{1},\bar{2}\rbrace.$

\section{ WHEN $\boldsymbol{\tau (R) \simeq C(R)}$?}

In order to show that $\tau (R)$ and $C(R)$ are isomorphic, it is
sufficient to find a bijection $f:R \longrightarrow R$ satisfying the following
condition:
$$a+b\in Z(R) \Longleftrightarrow f(a)-f(b)\in Z(R)$$
Obviously, if $2\notin Z(R)$, then $\tau (R) \not \simeq C(R)$
because by Lemma 2.4, $2\notin Z(R)$ implies that $\tau (R)$ is
not a regular graph, while $ C(R)$ is always regular. Hence $2\in
Z(R)$ is a necessary condition for the isomorphism $\tau (R)
\simeq C(R)$. Meanwhile, if $R$ is a ring with characteristic 2,
then $\tau (R) \simeq C(R)$; since $f={\mbox{id}}_{R}$ is the
desired graph isomorphism because $x=-x$ for every $x\in R$, and
$$a+b\in Z(R) \Longleftrightarrow f(a)-f(b)=a-b=a+b\in Z(R).$$
Moreover, if $R$ is a local ring (i.e., $Z(R)$ is an ideal of $R$
) and $2\in Z(R)$, then again $f={\mbox{id}}_{R}$ gives the
isomorphism $\tau (R) \simeq C(R)$; since $a+b \in Z(R)$ implies
$a+b-2b\in Z(R)$, which means that $f(a)-f(b)=a-b \in Z(R)$.

In order to generalize the two previous cases, let $R\cong R_{1}
\oplus \cdots \oplus R_{k}$ be a ring such that each $R_{i}$ is a
local ring of even order. Thus by Lemma 2.5, we have
$2_{R_{i}}\in Z(R_{i})$ for each $ i=1,\ldots ,k$. In this case,
if $a=(a_{1},\ldots ,a_{k})$ and $ b=(b_{1},\ldots ,b_{k})$ are
elements of $R$ and $a+b\in Z(R)$, then there exists an $i$ with
$1\le i\le k$ such that $ a_{i}+b_{i}\in Z(R_{i})$. Hence by the
same reasoning, we have $ a_{i}+b_{i}-2b_{i}=a_{i}-b_{i}\in
Z(R_{i})$, which implies that $ a-b\in Z(R)$. Therefore, if
$R\cong R_{1} \oplus \cdots \oplus R_{k}$ with $k\ge 1$, and every
$R_{i}$ is a local ring of even order, then $f={\mbox{id}}_{R}$
gives the isomorphism $\tau (R) \simeq C(R)$. The reader should
verify that this case also covers the first case when char$(R)=2$.

Now suppose that $R\cong R_{1} \oplus \cdots \oplus R_{k}$ such
that every $R_{i}$ is local and $F_{1}=R_{1}/\frak{M}_{1} \cong
\Bbb{Z}_{2}$. Then if $u_{1}$ and $ v_{1}$ are units in $R_{1}$,
we have $u_{1}+v_{1}\in \frak{M}_{1}=Z(R_{1}).$
 In this case, put
$$A=\lbrace (a_{1},\ldots ,a_{k})\in R\;|\;a_{1}\in Z(R_{1})\rbrace$$
and define $f:R \longrightarrow R$ by
$$f(x)=\bigg\{
\begin{gathered}
x\hspace{10mm} x\in A \\
-x \hspace{7mm}x\in A. \\
\end{gathered}
$$
Let $a=(a_{1},\ldots ,a_{k})$ and $b=(b_{1},\ldots ,b_{k})$ be
elements of $R$ such that $a+b\in Z(R)$. Then there are three
possibilities:
\begin{enumerate}
\item[{\rm(i)}]
 $a,b\in A$; then $f(a)-f(b)=a-b=(a_{1}-b_{1},\ldots
,a_{k}-b_{k}) $ belongs to $Z(R)$ since $Z(R_{1})$ is an ideal of
$R_{1}$.
\item[{\rm(ii)}] $a\in A$ and $b\notin A$ or vice versa; then $f(a)-f(b)=a+b$
which is already in $Z(R)$.
\item[{\rm(iii)}] $a,b\notin A$; then $a_{1},b_{1}\in R_{1}^{*}$, and since $
F_{1}=R_{1}/\frak{M}_{1} \cong \Bbb{Z}_{2}$, we have
$-a_{1}+b_{1}\in Z(R_{1})$. Hence in this case, we have
$f(a)-f(b)=-a+b$ which is also in $Z(R).$

\end{enumerate}
\noindent
 Consequently, $f$ gives the desired isomorphism $\tau
(R)\simeq C(R)$, since each possibility implies $f(a)-f(b)\in
Z(R)$ and $f$ is a bijection.

However, if $R\cong R_{1}\oplus \cdots \oplus R_{k}$, $2\in Z(R)$,
but $ F_{i}=\frac{R_{i}}{\frak{M}_{i}} \not \cong \Bbb{Z}_{2}$
for each $i=1,\ldots ,k$, then there is no guarantee for the
existence of an isomorphism $\tau (R) \simeq C(R)$.

\begin{figure}[h]
\centering
\includegraphics[width=6.7cm,height=6.02cm]{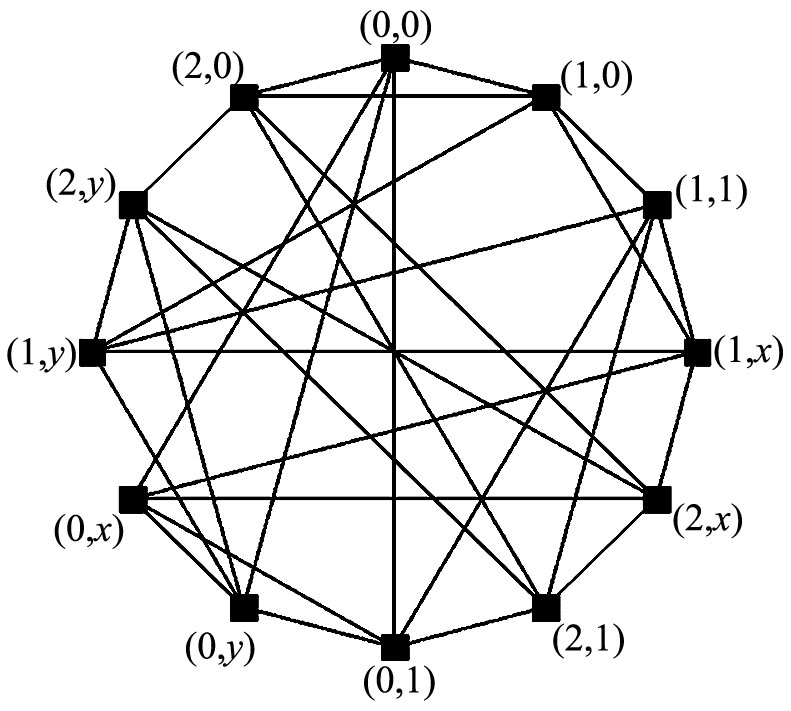}
\end{figure}
\begin{center}\textbf{Figure 3.} $C(\Bbb{Z}_{3}\oplus \Bbb{F}_{4})$.\end{center}
\begin{figure}[h]
\centering
\includegraphics[width=13.34cm,height=6.4cm]{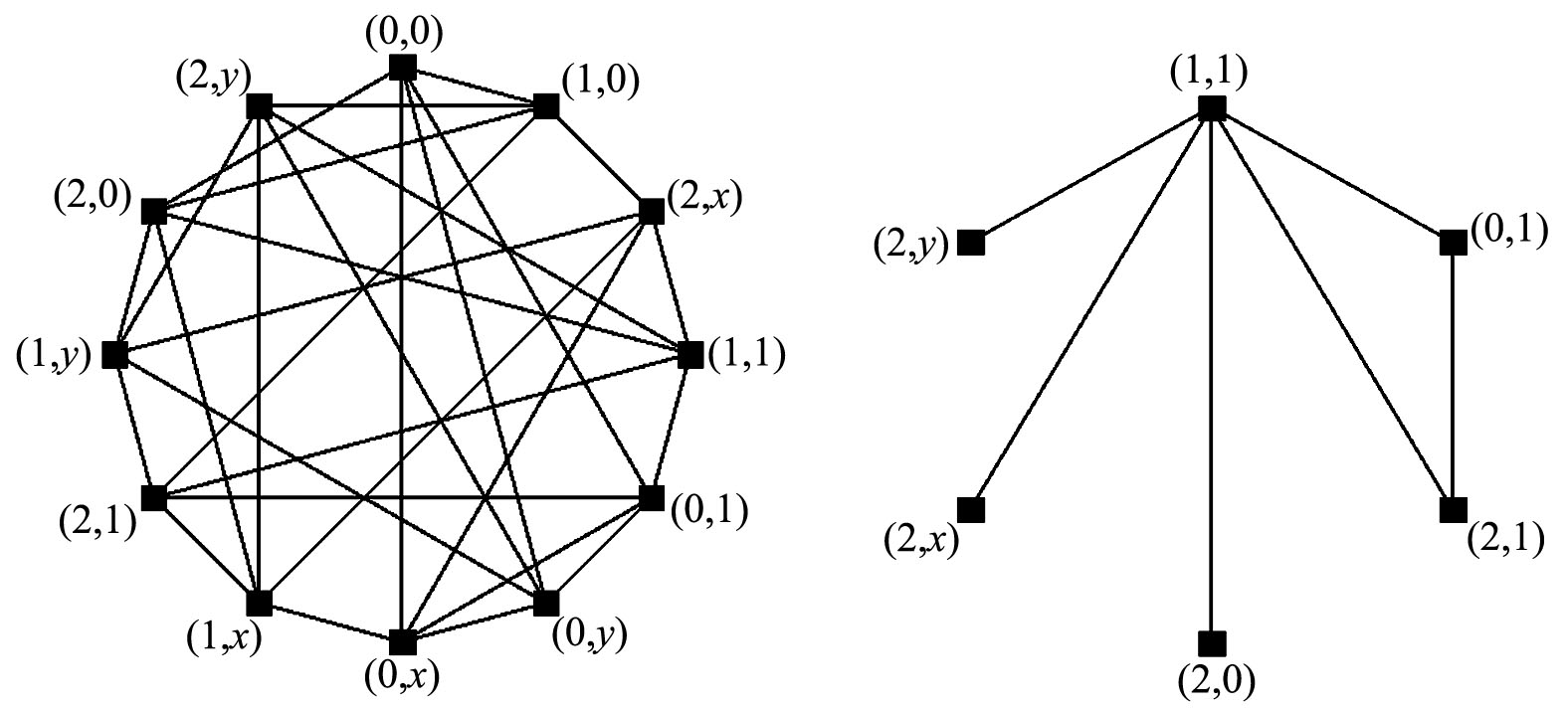}
\end{figure}
\begin{center}\textbf{Figure 4.} $\tau (\Bbb{Z}_{3}\oplus\Bbb{F}_{4})$ and its induced
subgraph $\tau (\Bbb{Z}_{3}\oplus \Bbb{F}_{4})\lbrack N\lbrack
(1,1)\rbrack \rbrack $.
\end{center}
\noindent
\textbf{Example 5.1.} $\tau (\Bbb{Z}_{3} \oplus
\Bbb{F}_{4})\not \simeq C(\Bbb{Z}_{3} \oplus \Bbb{F}_{4})$. To
see this,  suppose $\Bbb{F}_{4}=\{0,1,x,y\}$. Then the graph
$C(\Bbb{Z}_{3} \oplus \Bbb{F}_{4})$ has three disjoint 4-cliques
(see Figure 3):
\begin{eqnarray*}
 &\lbrace  (0,0),(0,1),(0,x),(0,y)\rbrace,& \\
 &\lbrace (1,0),(1,1),(1,x),(1,y)\rbrace,& \mbox{ and} \\
 &\lbrace (2,0),(2,1),(2,x),(2,y)\rbrace.&
\end{eqnarray*}
Hence every vertex in $C(\Bbb{Z}_{3}\oplus \Bbb{F}_{4})$ belongs
to a 4-clique. But the vertex $(1,1)$ in $\tau (\Bbb{Z}_{3} \oplus
\Bbb{F}_{4})$ is not a vertex of a 4-clique (see Figure 4).
Therefore, $\tau (\Bbb{Z}_{3} \oplus \Bbb{F}_{4})\not \simeq
C(\Bbb{Z}_{3}\oplus \Bbb{F}_{4})$.
\vspace{2mm}\\
\noindent
 \textbf{Theorem 5.2.} \textit{Let $R$ be a finite
commutative ring. Then the two graphs $\tau (R)$ and $C(R)$ are
isomorphic if and only if at least one of the following
conditions is true:
\begin{enumerate}
\item[{\rm(a)}]  $R\cong R_{1}\oplus \cdots \oplus R_{k}, k\ge 1$, and
each $R_{i}$ is a local ring of an even order,
\item[{\rm(b)}]  $R\cong
R_{1}\oplus \cdots \oplus R_{k}, k\ge 2$, and each $R_{i}$ is a
local ring and $f_{1}=2$.
\end{enumerate}}
\vspace{2mm} \noindent \textbf{Proof.} By the above discussion, if
at least one of (a) or (b) is true, then $\tau (R) \simeq C(R)$.
So, we only prove the converse. Suppose (a) and (b) do not hold
for a ring $R$. If $2\notin Z(R)$, then by the discussion at the
beginning of this section, we have $\tau (R) \not \simeq C(R).$
So, suppose that $2\in Z(R)$. Therefore, $R\cong R_{1}\oplus
\cdots \oplus R_{k}$; $k\ge 2$ or else $R$ is local and $2\in
Z(R)$ implies that condition (a) is true. Furthermore, we have $
R_{i}/Z(R_{i})\not \cong \Bbb{Z}_{2}$ for all $i=1,\ldots ,k$,
but since $2\in Z(R)$ , it can be inferred that
$R_{i}/Z(R_{i})\cong \Bbb{F}_{2^{t}}$, for some $ i=1,\ldots ,k$
and $t\ge 2$. Moreover, there exists an $i$, $ 1\le i\le k$, such
that $\vert R_{i}\vert $ is odd, because otherwise condition (a)
is true again. So, we assume that $\vert R_{1}\vert ,\ldots
,\vert R_{j}\vert $ are all even and $\vert R_{j+1}\vert ,\ldots
,\vert R_{k}\vert $ are all odd, after appropriate permutation of
factors of course.

In order to prove that $\tau (R)$ and $C(R)$ are not isomorphic,
we use a method similar to the proof of Example 5.1. In fact, to
prove that these two graphs are not isomorphic, we consider
maximal cliques that contain an edge.

If $R_{i}/Z(R_{i})\cong \Bbb{F}_{2^{t}}$ for some $t\ge 2$, then
every vertex $a=(a_{1},\ldots ,a_{k})$ belongs to a clique of
maximal size $ \vert R\vert /f_{i}$ in both graphs $\tau (R)$ and
$C(R)$ because $ \pi _{i}(a_{i})=-\pi _{i}(a_{i})$. But besides
belonging to these distinct maximal cliques, for each
$l=j+1,\ldots ,k$, the vertex $a$ belongs to a maximal $(|R|/f_{l}
) - $clique in $C(R)$. Because for each $ b=(b_{1},\ldots
,b_{k})$ and $c=(c_{1},\ldots ,c_{k})$ which are adjacent to $a$,
provided that $a_{l}-b_{l}\in Z(R_{l})$ and $ a_{l}-c_{l}\in
Z(R_{l})$, we have also $b_{l}-c_{l}\in Z(R_{l})$. This means
that $b$ and $c$ are also adjacent in $C(R)$. Therefore, if $b$
is adjacent to $a$, depending on the index $i$ that $
a_{i}-b_{i}\in Z(R_{i})$, we can say that $a$ and $b$ are both
vertices of a maximal $(\vert R\vert /f_{i})-$clique in $C(R)$.
Therefore, every edge of $C(R)$ belongs to a maximal $(\vert
R\vert /f_{i})-$clique, for some $i=1,\ldots ,k$.

Now, in $\tau (R)$, put $x=(0_{R_{1}},\ldots
,0_{R_{j}},-1_{R_{j+1}},\ldots ,-1_{R_{k}})$. Then $x$ is adjacent
to $1=(1,\ldots ,1)$, and for all $i=1,\ldots ,j$, we have $
x_{i}+1\notin Z(R_{i})$ but $x_{i}+1\in Z(R_{i} )$ for $i=j+1,\ldots ,k$.
We show that for $i=1,\ldots ,k$, the edge $\{1,x\}$ does not belong
to a maximal $(\vert R\vert /f_{i})-$clique in $\tau (R)$ and this
completes the proof.

Let $\{y_{s}\vert s \in S\}$ be a set of elements of $R$ of
maximal size which are adjacent to both $1$ and $x$ and also to
themselves. And, let $[a_{m}]$ denote the equivalence class of
$Z(R_{m})+a_{m}$ . If $\{y_{s}\vert s\in S\} \cup \{x,1\}$ forms a
clique of maximal size $|R|/f_{i}$, then there must be $1\le
m_{1}<m_{2}<\ldots <m_{q}\le k$; $0\le q\le k$ such that all
$y_{s}$'s belong to
$$R_{1}\oplus \cdots \oplus R_{m_{1}-1}\oplus \lbrack
a_{m_{1}}\rbrack \oplus R_{m_{1}+1}\oplus \cdots \oplus
R_{m_{q}-1}\oplus \lbrack a_{m_{q}}\rbrack \oplus R_{m_{q}+1} \oplus \cdots
\oplus R_{k}.$$
The case $q=0$ is not possible. To the contrary, suppose that $q=0$. Consequently, we must have

$$|R|/f_{i} =\bigg\{
\begin{gathered}
k\hspace{12mm} k \mbox{ is even} \\
k+1 \hspace{6mm} k \mbox{ is odd,} \\
\end{gathered}
 $$
for some $i=1,\ldots ,k$. But, if $k>2$, then the inequality $
k+1<3^{k-1}$ holds, and hence $|R|/f_{i}<3^{k-1}\le |R|/f_{i}$,
which is a contradiction. The case $k=2$ also implies a similar
contradiction.

If a vertex $y=(y_{1},\ldots ,y_{k})$ is adjacent to both $1$ and $
x$ in $\tau (R)$, then because each $f_{i}\ge 3$ and $\vert
R_{j+1}\vert ,\ldots ,\vert R_{k}\vert $ are all odd, there is not an $
i$, $1\le i\le k$, such that $y_{i}+1$ and $y_{i}+x_{i}\in Z(R_{i})$.
Hence if $y=(y_{1},\ldots ,y_{k})$ is adjacent to both $1$ and $x$,
 then there must exist $s,t \in \lbrace 1,\ldots ,k\rbrace $ with $s\ne t
$ such that $y_{s}+1\in Z(R_{s})$ and $y_{t}+x_{t}\in Z(R_{t})$.
Therefore, $q\ne 1$.

Now, suppose that $2\le q\le k$. We must have $\lbrack
a_{m_{p}}\rbrack =\lbrack -1_{m_{p}}\rbrack $ and $\lbrack
a_{m_{v}}\rbrack =\lbrack -x_{m_{v}}\rbrack $ for some $v\ne p$, $
1\le p,v\le k$, because otherwise the $y_{s}$'s are not adjacent
to $1$ or $x$. Meanwhile, for some $t\in \{1,\ldots ,q\},$ we must
also have $m_{t}\le j$, because otherwise the $y_{s}$'s cannot be
adjacent to themselves. Then, for each $i=1,\ldots ,k$ we have
$\vert \{y_{s} \vert s\in S\} \vert +2=\frac{\vert R\vert
}{\prod_{i=1}^{q}{f_{m_{i}}}}+2\ne \frac{\vert R\vert }{f_{i}}$.
To the contrary, suppose that $\frac{\vert R\vert
}{\prod_{i=1}^{q}{f_{m_{i}}}}+2=\frac{\vert R\vert }{f_{i}}$ for
some $ i=1,\ldots ,k$. Let $n$ be the greatest integer such that $
2^{n}{\big \vert} \vert R\vert $, i.e., $\vert R\vert =2^{n}r$ for
some odd integer $r\ne 1$. Since $R_{m_{1}},\ldots ,R_{m_{t}}$ are
local rings of even orders, it can be inferred that $
\prod_{i=1}^{t}{f_{m_{i}}}=2^{g}$ for some positive integer $g$, $
2\le g\le n$. Thus, we have
$$
\frac{2^{n}r}{2^{g}\prod_{i=t+1}^{q}{f_{m_{i}}}}+2=\frac{2^{n-g}r}
{\prod_{i=t+1}^{q}{f_{m_{i}}}}+2=2(\frac{2^{n-g-1}r}{\prod_{i=t+1}^{q}{f_{m_{i}}}}+1)=\frac{2^{n}r}{f_{i}}.$$
Hence,
$$
\frac{2^{n-g-1}r}{\prod_{i=t+1}^{q}{f_{m_{i}}}}+1=\frac{2^{n-1}r}{f_{i}}
.$$

If $n-g-1>0$, then the previous equality is impossible since the
right hand side is always an even integer but the left hand side is odd.
The case $n-g-1=0$ is not possible, because it means that one of the $
f_{i}$'s equals $2$, a contradiction to our assumption. The case $
n-g-1=-1$ cannot satisfy the previous equality, because the left hand
side is not an integer, while the right hand side is always an integer.
Consequently, the case $2\le q\le k$ is not possible either.

Therefore, in $\tau (R)$, the set $\{y_{s}\vert s\in S\} \cup
\{x,1\} $ cannot form a clique of maximal size $\frac{\vert
R\vert}{f_{i}}$ for each $i$, $i=1,\ldots ,k$. Hence $\tau (R)\not
\simeq C(R)$. $\Box$
\vspace{2mm}\\
\textbf{Remark 5.3.} The two conditions (a) and (b) in Theorem 5.2
are not distinct, i.e., they can both be true for a ring $R$.
\section*{ACKNOWLEDGMENTS}

The authors owe a great debt to the referees who have carefully
read an earlier version of this paper and made significant
suggestions for improvement. We would like to express our deep
appreciation for the referees' work.

\section*{REFERENCES}

Akbari, S., Kiani, D., Mahammadi, F., Moradi, S. (2009). The total graph
and regular graph of a commutative ring, \textit{J. Pure Appl. Algebra
} 213: 2224-2228.
\vspace{2mm}\\
Akhtar, R., Bogges, M., Jackson-Henderson, T., Jiménez, I., Karpman, R.,
Kinzel, A., Pritikin, D. (2009). On the unitary Cayley graph of a finite
ring, \textit{Electronic J. of Combinatorics} 16: \#R117.
\vspace{2mm}\\
Anderson, D. F., Badawi, A. (2008). The total graph of a commutative
ring, \textit{J. Algebra} 217: 2706-2719.
\vspace{2mm}\\
Atiyah, M. F., Macdonald, I. G. (1969). Introduction to commutative
algebra, Addison-Wesley Publishing Co.
\vspace{2mm}\\
Behzad, M. (1970). A Characterization of total graphs, \textit{
Proceedings Amer. Math. Society} 26(3): 383-389.
\vspace{2mm}\\
Beineke, L. W., Wilson, R. J. (2004). Topics in algebraic graph theory,
Cambridge University Press.
\vspace{2mm}\\
Maimani, H. R., Wickham, C., Yassemi, S.  Rings whose total
graphs have genus at most one, \textit{Rocky Mountain J. Math.}:
to appear.
\vspace{2mm}\\
West, D. B. (2000). Introduction to graph theory, second edition,
Prentice-Hall.

\end{document}